\documentclass[12pt, a4paper, twoside]{amsart}
\usepackage{amssymb}
\usepackage{graphics}

\usepackage[utf8]{inputenc}
\title{The Hochschild cohomology of the group~$G^2_3$}

\author{Hassan AlHussein}
\address{Novosibirsk State University, Russia}
\author{Pavel Kolesnikov}
\address{Sobolev Institute of Mathematics, Russia}

\newtheorem{theorem}{Theorem}
\newtheorem{proposition}{Proposition}
\newtheorem{definition}{Definition}
\newtheorem{example}{Example}
\newtheorem{corollary}{Corollary}

\def\GSB{\mathrm{GSB}} 
\def\Ker{\mathop{\fam 0 Ker}\nolimits}
\def\rank{\mathop{\fam 0 rank}\nolimits} 
\def\Hom{\mathop{\fam 0 Hom}\nolimits} 

\begin{document}

\begin{abstract}
Cohomology plays an important role in algebra. The group of Hochschild
cohomology for group $G^2_3$ contains information about its structure. In this work we found group of Hochschild $n$-cohomology for the group $G^2_3$ for any $n\ge1$. Anick resolution, Algebraic Morse theory and Gr\"obner---Shirshov basis were used to get this result.   
\end{abstract}
\keywords{Hochschild cohomology, Anick resolution, Gr\"obner--Shirshov basis, Morse matching}

\maketitle

\section{Introduction}
Homological methods allow us to get important information about the structure of an algebra. For associative algebras, Hochschild cohomologies play an important role in structure and representation theory. Finding the
Hochschild cohomology group $H^n(A,M)$ of a given algebra $A$ with coefficients in a given $A$-bimodule $M$ is often a difficult problem.
In order to solve this problem one needs a long exact sequence 
starting from $A$, a resolution of $A$.
The most natural bar-resolution is easy to construct but it is too bulky 
for computations. Another approach was proposed by David J.~Anick in 1986 \cite{Anick}, where it was built a free resolution for associative algebra
which is homotopy equivalent to the bar-resolution. The Anick resolution 
was also used to find Poincare Series \cite{ufn}.
Computation of the differentials in the Anick resolution according to 
the original algorithm described in \cite{Anick} is extremely hard. 
In order to make the computation easier, one may use the discrete algebraic
Morse theory based on the concept of a Morse matching defined in \cite{jollwel}. This concept was used in geometry first, then it became applicable in algebra \cite{formancell, formanguide}. 

A series of groups $G_n^k$ were introduced in \cite{manturov,manturov1,manturov2,manturov3} in the study of braids, links, and Coxeter groups. 
In this paper, we consider the group $G_3^2$ originated in the following 
(informal) context. Assume three points are moving without collisions on a disk in the plane at time $t\in [0,1]$. The trajectories of these three points are characterized by an element of a group constructed as follows. Whenever the points stand in one line, say, the point $i$ stands between 
$j$ and $k$ ($\{i,j,k\}=\{1,2,3\}$), we write down the generator $a_{jk}$.
Therefore, we have three generators 
$a=a_{12}$, $b=a_{13}$, and $c=a_{23}$, and the trajectories of three points are characterized by a word in the alphabet $\{a,b,c\}$. 
The following relations hold on these generators:
\[
a^2=b^2=c^2=1, \quad bac=cab.
\]

In the present work, we apply the Morse matching theory to find the Anick resolution and 
calculate the groups of Hochschild $n$-cohomologies of the group $G^2_3$ for all $n\ge 1$ with coefficients in all 1-dimensional 
$kG^2_3$-bimodules. It is easy to see that all such modules are of the form $W={}_{\epsilon_1} W_{\epsilon_2}$, where 
\[
{}_{\epsilon_1}W_{\epsilon_2} = \Bbbk 1, 
\quad x\cdot 1 = \epsilon_1(x),\quad 1\cdot x = \epsilon_2(x)
\]
for $x\in \{a,b,c\}$, $\epsilon_i: \{a,b,c\} \to \{1,-1\}$.

\section{Morse matching}

In this section, we state main definitions and essential results related to the construction of the Anick chain complex via algebraic Morse theory following \cite{Anick, lopatkin, Maclane, Bokut, skold, jollwel}. Let $k$ be
a field and let $\Lambda$ be a unital associative  $k$-algebra with  an augmentation, i.e., a $k$-algebra homomorphism $\varepsilon:\Lambda \rightarrow k$. Let $A$ be a set of generators for $\Lambda$. Suppose that  $\le$ is a well ordering on the free monoid $A^{*}$ generated by~$A$. 
Denote by $k \langle A \rangle$ the free associative algebra with identity generated by~$A$. There is a canonical surjection 
$f:k \langle A \rangle \to \Lambda$, so that 
$\Lambda \simeq k \langle A \rangle / \Ker(f)$.
Let $\GSB_\Lambda$ be a Gr\"obner---Shirshov basis of $\Ker f$,
i.e., a confluent set of defining relations for the algebra $\Lambda $.
Denote by $V$ the set of the leading terms of relations from $\GSB_\Lambda$. Following Anick \cite{Anick}, $V$ is called the set of {\em obstructions}.

For $n \ge 1$, a word $v=x_{i_1}\ldots x_{i_t} \in A^{*}$ 
is an $n$-{\em prechain} if and only if there exist 
$a_j,b_j \in \mathbb{Z}$, $1 \le j \le n$, satisfying the following conditions:
\begin{enumerate}
    \item $1=a_1<a_2 \le b_2<\ldots<a_n \le b_{n-1}<b_n=t$;
    \item $ x_{i_{a_j}}\ldots x_{i_{b_j}} \in V$ for $1 \le j \le n$.
\end{enumerate}
An $n$-prechain $x_{i_1}\ldots x_{i_t}$ is an $n$-{\em chain} 
if only if the integers $a_j,b_j$ can be chosen in such a way that  $x_{i_1}\ldots x_{i_t}$  is not an $m$-prechain for neither 
$s < b_m$, $1 \le m \le n$.

As in \cite{Anick}, we say that the elements of $A$ are $0$-chains, the elements of $V$ are $1$-chains, and denote the set of $n$-chains by $V^{(n)}$.
The cokernel of  $\varepsilon: k \to \Lambda$ is denoted by $\Lambda\slash k$.
A word $w\in A^*$ is said to be $V$-{\em reduced} if it does not contain a 
word from $V$ as a subword. Since $V$ is the set of all obstructions, the  set of all non-trivial $V$-{\em reduced} words forms a linear basis of $\Lambda/k$
\cite{Bok72}.

The two-sided bar resolution $(B_n(\Lambda,\Lambda),d_n)_{n\ge 0}$ is a $\Lambda$-free bimodule resolution of $\Lambda$, where
$$
B_n(\Lambda,\Lambda):=\Lambda\otimes (\Lambda\slash k )^{\otimes n}\otimes \Lambda.
$$
An element $1\otimes \lambda_1\otimes \dots \lambda_n\otimes 1\in B_n(\Lambda,\Lambda )$
is denoted by $[\lambda_1|\ldots |\lambda_n]$ for $\lambda_i\in \Lambda/k$.
The differential $d_n$ is defined as follows:
\begin{multline}\nonumber
 d_n([\lambda_1|\ldots|\lambda_n])=\lambda_1[\lambda_2|\ldots|\lambda_n]+(-1)^n[\lambda_1|\ldots|\lambda_{n-1}]\lambda_n+
 \\
 +\sum\limits_{i=1}^{n-1}(-1)^i[\lambda_1|\ldots|\lambda_i\lambda_{i+1}|
 \ldots|\lambda_n]
\end{multline} 

Let $B=(B_n,d_n)_{n \ge 0}$  be a chain complex of free (left) $\Lambda$-modules.
Choose a basis $X_n$ in each $B_n$, and set $X=\cup_{n>0}X_n$. Then $B_n={\oplus_{b \in X_n}}\Lambda b$, and the differential $d_n:B_n\to B_{n-1}$ may be uniquely presented in the following form:
 $$
 d_n(b)=\sum\limits_{b' \in X_{n-1}}[b:b']b',\quad b\in X_n,
    $$
where $[b:b']\in \Lambda $.
Construct a directed weighted graph $\Gamma(B)=(X,E)$
considering $X$ as the set of vertices.
Edges of $\Gamma(B)$ 
are ordered triples of the form 
$e=(x,x',\lambda )\in X\times X\times \Lambda $, where $x$ is the source of $e$, 
$x'$ is the range of $e$, and $\lambda =\omega(e)$  is the weight of~$e$.
Define the set $E$ of weighted edges by the following rule:
$$
 (b,b',[b,b']) \in E \iff b \in X_i,\ b' \in X_{i-1},\ [b:b'] \ne 0.
$$
 
\begin{definition}\label{defn1.1}
A finite subset $M \subset E$ is called a {\em Morse matching} if and only if
\begin{itemize}
\item Each vertex $x\in X$ lies in at most one edge $e\in M$;
\item For all edges $(b,b',[b,b']) \in M$, the weight $[b,b']$
is an invertible element in the center of $\Lambda $;
\item The graph $\Gamma_M(X,E_M)$ has no directed cycles, where $E_M$ is given by
$$
E_M=(E \setminus M)\cup\{(b',b,-[b:b']^{-1})\mid (b,b',[b:b'])\in M\}.
$$
\end{itemize}
\end{definition}

A vertex $b\in X$ is {\em critical} with respect to $M$ if $b$ does not 
belong to an edge $e\in M$; we denote by $X_n^M$ the set of critical edges from $X_n$. 

Suppose $p$ is a path $p=b_1 \rightarrow \ldots \rightarrow b_r$ in a 
weighted directed graph with vertices $X$. Then 
$$
  \omega(p):=\prod\limits_{i=1}^{r-1}\lambda_i,
 \quad
  e_i = (b_i,b_{i+1},\lambda_i).
$$
Denote by
$\Gamma(b,b')$, $b,b'\in X$, 
the sum of weights of all path from $b$ to~$b'$.

\begin{theorem}\label{homotopy}
A chain complex  $(B_n,d_n)_{n \ge 0}$  is homotopy-equivalent
to the complex $(B_n^M,d_n^M)_{n \ge 0}$ of free $\Lambda$-modules  where $B_n^M={\oplus_{b \in X_n^M}}\Lambda b$ and
\begin{equation}\label{eq:d_M-product}
d_n^M: B_n^M \rightarrow B_{n-1}^M; \quad 
d_n^M(b)=\sum\limits_{b' \in X_{n-1}^M}\Gamma(b,b')b',
\end{equation}
where $\Gamma(b,b')$ is calculated in the graph $\Gamma_M(B)=\Gamma_M(X,E_M)$.
\end{theorem}

Let $B=(B_n,d_n)_{n\ge 0}$ be the bimodule bar complex for $\Lambda $.
We may consider $B$ as a complex of left $\Lambda\otimes \Lambda^{op}$-modules and construct a Morse matching in $\Gamma (B)$
in the following way.

\begin{theorem}\label{thm:Matching}
For $w\in A^{*}$ let $V_{w,i}$, be the set of vertices $[w_1|\ldots|w_n]$ 
in $\Gamma(B)$ such that $w=w_1\ldots w_{n}$, all $w_k$ are $V$-reduced, 
and $i\ge-1$ is the largest integer  such that $w_1\ldots w_{i+1} \in V^{(i)}$.
Let $V_w ={\cup_{i\ge-1}}V_{w,i}$; define a partial matching $M_w$ 
by letting $M_w$ consist of all edges
$$
 [w_1|\ldots|w'_{i+2}|w''_{i+2}|\ldots|w_n] \rightarrow [w_1|\ldots|w_{i+2}|\ldots|w_n],
 $$
for 
$[w_1|\ldots|w_n] \in V_{w,i}$, where
$w'_{i+2}w''_{i+2}=w_{i+2}$ and $[w_1|\ldots|w_{i+1}|w'_{i+2}] \in V^{(i+1)}$. 
Then the set of edges $ M={\cup_{w}}M_w$ is a Morse matching on $\Gamma(B)$ with critical cells $X_n^M=V^{(n-1)}$ for all~$n$.
\end{theorem}

\begin{proposition}\label{prop:AnickComplex} 
The chain complex $(B_{n}^M(\Lambda),d^M_n)_{n\ge 0}$
is the $\Lambda$-free Anick resolution.
\end{proposition}

\begin{definition}
Let $\Lambda $ be a unital $k$-algebra, and let $W$ be a $\Lambda $-bimodule. 
Consider the cochain complex
$$
 0 \rightarrow M \stackrel{\Delta_0}{\rightarrow} \Hom_k(\Lambda, W)
 \stackrel{\Delta_1}{\rightarrow} 
 \ldots \stackrel{\Delta_{n-1}}{\rightarrow} \Hom_k(\Lambda^{{\otimes n}},W) \stackrel{\Delta_n}{\rightarrow} \ldots
 $$
where
\begin{multline}\nonumber
\Delta_n f(r_1, \ldots, r_{n+1})=r_1 f(r_2, \ldots, r_{n+1}) \\
+ \sum\limits_{i=1}^{n}(-1)^i f(r_1,\ldots,r_i r_{i+1},\ldots,r_{n+1}) +
 (-1)^{n+1}f(r_1,\ldots,r_n)r_{n+1},
\end{multline}
$f \in \Hom_k(\Lambda ^{{\otimes n}},W)$.
\end{definition}
The $n$th Hochschild cohomology group of $\Lambda $ with coefficients in $W$ is
$H^n(\Lambda ,W)={\Ker \Delta_n}/{\mathrm{Im}\,\Delta_{n-1}}$.

We may use the Anick resolution to calculate $H^n(\Lambda ,W)$ by means of the following diagram:
$$
 \begin{picture}(60,60)
 \put(-80,50)
 {$
  \ldots \longleftarrow B_{n-1} \stackrel{d_{n}} \longleftarrow B_n
  \stackrel{d_{n+1}} \longleftarrow B_{n+1} \longleftarrow \ldots
  $}
 \put(-33,42){\vector(0,-1){30}}
 \put(-45,27){{$f$}}
 \put(-40,0){$W$}
 \put(10,42){\vector(-1,-1){30}}
 \put(-22,30){\footnotesize{$\Delta_{n}f$}}
 \put(20,42){\vector(0,-1){30}}
 \put(8,27){$g$}
 \put(12,0){$W$}
 \put(55,42){\vector(-1,-1){30}}
 \put(45,25){\footnotesize{$\Delta_{n+1}g$}}
 \end{picture}
 $$
 \noindent
Here 
\[
\begin{gathered}
\beta^n:=\mathrm{Im}\,\Delta_{n}=\{g:B_n \rightarrow W 
  \mid g=\Delta_{n}f=f d_{n}; \, f: B_{n-1} \rightarrow W\}, \\ 
Z^n:=\Ker\Delta_{n+1}=\{g:B_n \rightarrow W; \, \Delta_{n+1}g=g d_{n+1}=0\}, \end{gathered}
\]
so $H^n(\Lambda,W)=Z^n \slash \beta^n$.

\section{The Anick complex for $G^2_3$}

\begin{definition}[\cite{manturov}]\label{defn1.3}
The group $G^2_3$ is defined by generators $a, b, c$, and relations $S_2$, i.e., $G^2_3=\mathrm{Sgp}\,\langle a, b, c \mid S_2\rangle $,
where
\[
S_2 = \{a^2-1, b^2-1, c^2-1, bca-acb \}.
\]
\end{definition}

Obviously, the set $S_2$ defines the same ideal in 
$k\langle a,b,c\rangle $ as $S=\{a^2-1, b^2-1, c^2-1, ba-cabc, bca-acb\}$.

\begin{theorem}[\cite{YUXIU}]\label{thm:GSB}
The set of relations $S$ is a confluent set of rewriting rules.
\end{theorem}

To be more precise, Theorem \ref{thm:GSB} states that $S$ is a 
Gr\"obner---Shirshov basis of the group algebra $kG^2_3$
with respect to the tower order on 
$\{a, b, c \}^*$, where $a>b>c$.

To find Anick complex for $G^2_3$, we need several steps.
First, we have to find the set of obstructions for the group $G^2_3$ 
relative to the given Gr\"obner---Shirshov basis (the set of leading terms in $S=\GSB_{G^2_3}$) 
and the set of $n$-chains for $n\ge 0$. 
Next, we should build a Morse matching in the graph $\Gamma(B)$
and construct the graph $\Gamma_M(B)$. 
Finally, it remains to calculate the Anick differential.

In order to make notations shorter, we will often use 
$[w]$ for an element $[w_1|\dots |w_n]$, where $w=w_1\dots w_n$
is an ($n-1$)-chain.

\begin{example}\label{exmp:d_1}
For $n=1$, 
$V^{(1)}=\{a^2, b^2, c^2, ba, bca\}$. Let us construct Morse matching
 and evaluate $d_2: B_2\to B_1$ as described in Theorem \ref{thm:Matching}.
The simplest case is 
$$
   [a] \stackrel{a\otimes 1} \longleftarrow [a|a] \stackrel{1\otimes a} \longrightarrow [a]
$$
Hence, $d_2[a^2]=a[a]+[a]a$. Similarly, 
$d_2[b^2]=b[b]+[b]b$, $d_2[c^2]=c[c]+[c]c$.
A more complicated construction is needed for $d_2[b|ca]$:

\newpage

\centerline{\hspace{150mm}\vbox{%
\begin{picture}(60,-100)
  \put(-80,-25)
  {$
   [ca] \stackrel{b\otimes 1} \longleftarrow [b|ca] \stackrel{1\otimes ca} \longrightarrow [b]
   $}
  \put(-85,-22){\vector(-1,0){20}}
  \put(-105,-26){\vector(1,0){20}}
  \put(-100,-35){\footnotesize{$-1$}}
  \put(-95,-20){\footnotesize{$1$}}
  \put(-130,-25){$[c|a]$}
  \put(-120,-30){\vector(0,-1){30}}
  \put(-115,-50){\footnotesize{$1\otimes a$}}
  \put(-125,-70){$[c]$}
  \put(-130,-22){\vector(-1,0){40}}
  \put(-155,-20){\footnotesize{$c\otimes 1$}}
  \put(-183,-25){$[a]$}
  \put(-25,-30){\vector(0,-1){20}}
  \put(-42,-40){\footnotesize{$-1$}}
  \put(-36,-60){$[acb]$}
  \put(-28,-65){\vector(0,-1){20}}
  \put(-23,-84){\vector(0,1){20}}
  \put(-22,-77){\footnotesize{$-1$}}
  \put(-38,-77){\footnotesize{$1$}}
  \put(-81,-95)
  {$
   [cb] \stackrel{a\otimes 1} \longleftarrow [a|cb] \stackrel{1\otimes cb} \longrightarrow [a]
   $}
  \put(-75,-99){\vector(0,-1){20}}
  \put(-70,-118){\vector(0,1){20}}
  \put(-65,-110){\footnotesize{$-1$}}
  \put(-85,-110){\footnotesize{$1$}}
  \put(-120,-130)
  {$
   [b] \stackrel{c\otimes 1} \longleftarrow [c|b] \stackrel{1\otimes b} \longrightarrow [c]
   $}
\end{picture}}}

\vspace{140bp}

Hence,
$d_2[bca]=[b]ca+bc[a]+b[c]a-[a]cb-a[c]b-ac[b]$.
Similarly,
$d_2[ba]=b[a]+[b]a-[c]abc-c[a]bc-cab[c]-ca[b]c$.
\end{example}

\begin{example}\label{exmp:d_2}
For $n=2$, $V^{(2)}=\{a^3, b^3, c^3, b^2a, ba^2, b^2ca, bca^2\}$.
Let us construct Morse matching
and evaluate $d_3: B_3\to B_2$. The simplest case is
$$
   [a|a] \stackrel{a\otimes 1} \longleftarrow [a|a|a] \stackrel{-1\otimes a} \longrightarrow [a|a].
$$
Hence, $d_3[a^3]=a[a^2]-[a^2]a$. Similarly, 
$d_3[b^3]=b[b^2]-[b^2]b$, 
$d_3[c^3]=c[c^2]-[c^2]c$. 
A more complicated construction is needed for $d_3[bca^2]$:
\vspace{220bp}
\centerline{\hspace{150mm}\vbox{%
\begin{picture}(60,-100)
 \put(-60,-20)
  {$
  [ca|a]\stackrel{b\otimes 1}\longleftarrow [b|ca|a] \stackrel{-1\otimes a}\longrightarrow [b|ca]
  $}
  \put(10,-25){\vector(0,-1){20}}
  \put(-45,-25){\vector(0,-1){20}}
  \put(-50,-45){\vector(0,1){20}}
  \put(-42,-38){\footnotesize{$-1$}}
  \put(-58,-38){\footnotesize{$1$}}
  \put(-114,-55)
  {$
   [a|a]\stackrel{c\otimes 1}\longleftarrow [c|a|a] 
   $}
  \put(-5,-40){\footnotesize{$-1$}}
  \put(-8,-55){$[acb|a]$}
  \put(4,-60){\vector(0,-1){20}}
  \put(8,-80){\vector(0,1){20}}
  \put(8,-70){\footnotesize{$-1$}}
  \put(-8,-70){\footnotesize{$1$}}
  \put(-63,-90)
  {$
  [cb|a]\stackrel{a\otimes 1}\longleftarrow [a|cb|a]
  $}
  \put(10,-107){\vector(0,-1){20}}
  \put(12,-118){\footnotesize{$1$}}
  \put(-5,-140)
  {$
  [a|abc]
  $}
   \put(10,-145){\vector(0,-1){20}}
   \put(15,-165){\vector(0,1){20}}
   \put(18,-160){\footnotesize{$-1$}}
   \put(1,-160){\footnotesize{$1$}}
           \put(5,-185)
  {$
   [a|a|bc]\stackrel{-1\otimes bc}\longrightarrow[a|a] 
   $}
   \put(-49,-95){\vector(0,-1){20}}
   \put(-45,-115){\vector(0,1){20}}
   \put(-60,-107){\footnotesize{$1$}}
   \put(-40,-107){\footnotesize{$-1$}}
    \put(-114,-130)
  {$
         [b|a]\stackrel{c\otimes 1}\longleftarrow [c|b|a]
   $}
   \put(-50,-135){\vector(0,-1){20}}
   \put(-45,-147){\footnotesize{$1$}}
   \put(-72,-166)
  {$
   [c|cabc]
   $}
  \put(-54,-175){\vector(0,-1){20}}
  \put(-48,-195){\vector(0,1){20}}
  \put(-43,-185){\footnotesize{$1$}}
  \put(-70,-185){\footnotesize{$-1$}}
  \put(-70,-209)
  {$
   [c|c|abc]
   $}
  \put(-27,-208){\vector(1,0){35}}
  \put(-29,-205){\footnotesize{$-1\otimes abc$}}
  \put(10,-210)
  {$
   [c|c]
   $}
\end{picture}}}   
Hence, 
$d_3[bca^2]=bc[a^2]-[bca]a-ac[ba]-a[c^2]abc-[a^2]bc$.
Similarly,
$d_3[b^2a]=b[ba]-[b^2]a+[bca]bc+ac[b^2]c+a[c^2]$,
$d_3[ba^2]=b[a^2]-[ba]a-ca[bca]-c[a^2]cb-[c^2]b$,
$d_3[b^2ca]=b[bca]-[b^2]ca+[ba]cb+ca[b^2]+cab[c^2]b$.
\end{example}

\begin{theorem}\label{thm:Res}
For $n\ge3$, $V^{(n)}=\{a^{n+1}, b^{n+1}, c^{n+1}, b^{i}a^{j},  b^{i}ca^{j} \mid  i+j=n+1, i,j\ge1\}$.
Then the Anick differential $d_{n+1}: B_{n+1}\to B_n$ 
is given by the following rules.
\[
\begin{aligned}
d_{n+1}[a^{n+1}]&=a[a^n]+(-1)^{n+1}[a^n]a, \\
d_{n+1}[b^{n+1}]&=b[b^n]+(-1)^{n+1}[b^n]b, \\
d_{n+1}[c^{n+1}]&=c[c^n]+(-1)^{n+1}[c^n]c.
\end{aligned}
\]
When $n+1$ is an even number,
\[
\begin{aligned}
d_{n+1}[b^na]&=b[b^{n-1}a]+[b^n]a-[b^{n-1}ca]bc-ca[b^n]c-cab[c^n],\\
d_{n+1}[ba^n]&=b[a^n]+[ba^{n-1}]a-ca[bca^{n-1}]-c[a^n]bc-[c^n]abc,\\
d_{n+1}[b^nca]&=b[b^{n-1}ca]+[b^n]ca-[b^{n-1}a]cb-ac[b^n]-a[c^n]b,\\
d_{n+1}[bca^n]&=bc[a^n]+[bca^{n-1}]a-ac[ba^{n-1}]-[a^n]cb-a[c^n]b,
\end{aligned}
\]
if $j=2,4,\dots, n-1$ is even then 
\[
\begin{aligned}
d_{n+1}[b^ia^j] & =b[b^{i-1}a^j]+[b^ia^{j-1}]a+[b^{i-1}ca^j]cb+ac[b^ica^{j-1}],\\
d_{n+1}[b^ica^j]&=b[b^{i-1}ca^j]+[b^ica^{j-1}]a+[b^{i-1}a^j]bc+ca[b^ia^{j-1}]\\
&+\binom{\frac{n-3}{2}}{\frac{j-2}{2}}cab[c^n]abc-\binom{\frac{n-3}{2}}{\frac{i-2}{2}}[c^n],
\end{aligned}
\]
if $j=3,5,\dots, n-2$ is odd then
\[
\begin{aligned}
d_{n+1}[b^ia^j]&=b[b^{i-1}a^j]+[b^ia^{j-1}]a-[b^{i-1}ca^j]bc-ca[b^ica^{j-1}] \\
   &-\binom{\frac{n-3}{2}}{\frac{i-3}{2}}cab[c^n]-\binom{\frac{n-3}{2}}{\frac{j-3}{2}}[c^n]abc, \\
d_{n+1}[b^ica^j]&=b[b^{i-1}ca^j]+[b^ica^{j-1}]a-[b^{i-1}a^j]cb-ac[b^ia^{j-1}] \\
 & -\left\{\binom{\frac{n-3}{2}}{\frac{i-3}{2}}+\binom{\frac{n-3}{2}}{\frac{j-3}{2}}\right \}a[c^n]b.
\end{aligned}
\]
When ${n+1}$ is an odd number,
\[
\begin{aligned}
d_{n+1}[b^na]&=b[b^{n-1}a]-[b^n]a+[b^{n-1}ca]bc+ac[b^n]c+a[c^n], \\
d_{n+1}[ba^n]&=b[a^n]-[ba^{n-1}]a-ca[bca^{n-1}]-c[a^n]cb-[c^n]b, \\
d_{n+1}[b^nca]&=b[b^{n-1}ca]-[b^n]ca
   +[b^{n-1}a]cb+ca[b^n]+cab[c^n]b,\\
d_{n+1}[bca^n]&=bc[a^n]-[bca^{n-1}]a-ac[ba^{n-1}]-[a^n]bc
-a[c^n]abc,
\end{aligned}
\]
if $j=2,4,\dots, n-2$ is even  then
\[
\begin{aligned}
d_{n+1}[b^ia^j]&=b[b^{i-1}a^j]-[b^ia^{j-1}]a-[b^{i-1}ca^j]cb \\
& -ca[b^ica^{j-1}]-\binom{\frac{n-2}{2}}{\frac{j-2}{2}}[c^n]b,\\
d_{n+1}[b^ica^j]&=b[b^{i-1}ca^j]-[b^ica^{j-1}]a-[b^{i-1}a^j]bc \\ & -ac[b^ia^{j-1}]-\binom{\frac{n-2}{2}}{\frac{j-2}{2}}a[c^n]abc,
\end{aligned}
\]
if $j=3,5,\dots, n-1$ is odd then
\[
\begin{aligned}
d_{n+1}[b^ia^j]&=b[b^{i-1}a^j]-[b^ia^{j-1}]a+[b^{i-1}ca^j]bc\\
& +ac[b^ica^{j-1}]+\binom{\frac{n-2}{2}}{\frac{i-2}{2}}a[c^n],\\
d_{n+1}[b^ica^j]&=b[b^{i-1}ca^j]-[b^ica^{j-1}]a+[b^{i-1}a^j]cb \\
& +ca[b^ia^{j-1}]+\binom{\frac{n-2}{2}}{\frac{i-2}{2}}cab[c^n]b.
\end{aligned}
\]
\end{theorem}

\textbf{Proof}. It is enough to construct a Morse  matching $M$ for 
the graph $\Gamma (B)$ corresponding to the bar resolution and calculate the 
differentials $d_{n+1}=d_{n+1}^M$ in the graph $\Gamma_M(B)$. 
As in Examples \ref{exmp:d_1}, \ref{exmp:d_2}, this can be done 
by the algorithm described in Theorem \ref{thm:Matching}. Namely, for each $w\in V^{(n)}$ 
we need to construct the subgraph of $\Gamma_M(B)$ which contains all paths from $[w]$ to 
vertices $[v]$, $v\in V^{(n-1)}$, and calculate the differentials~\eqref{eq:d_M-product}.

\begin{figure}\label{fig1}
\includegraphics{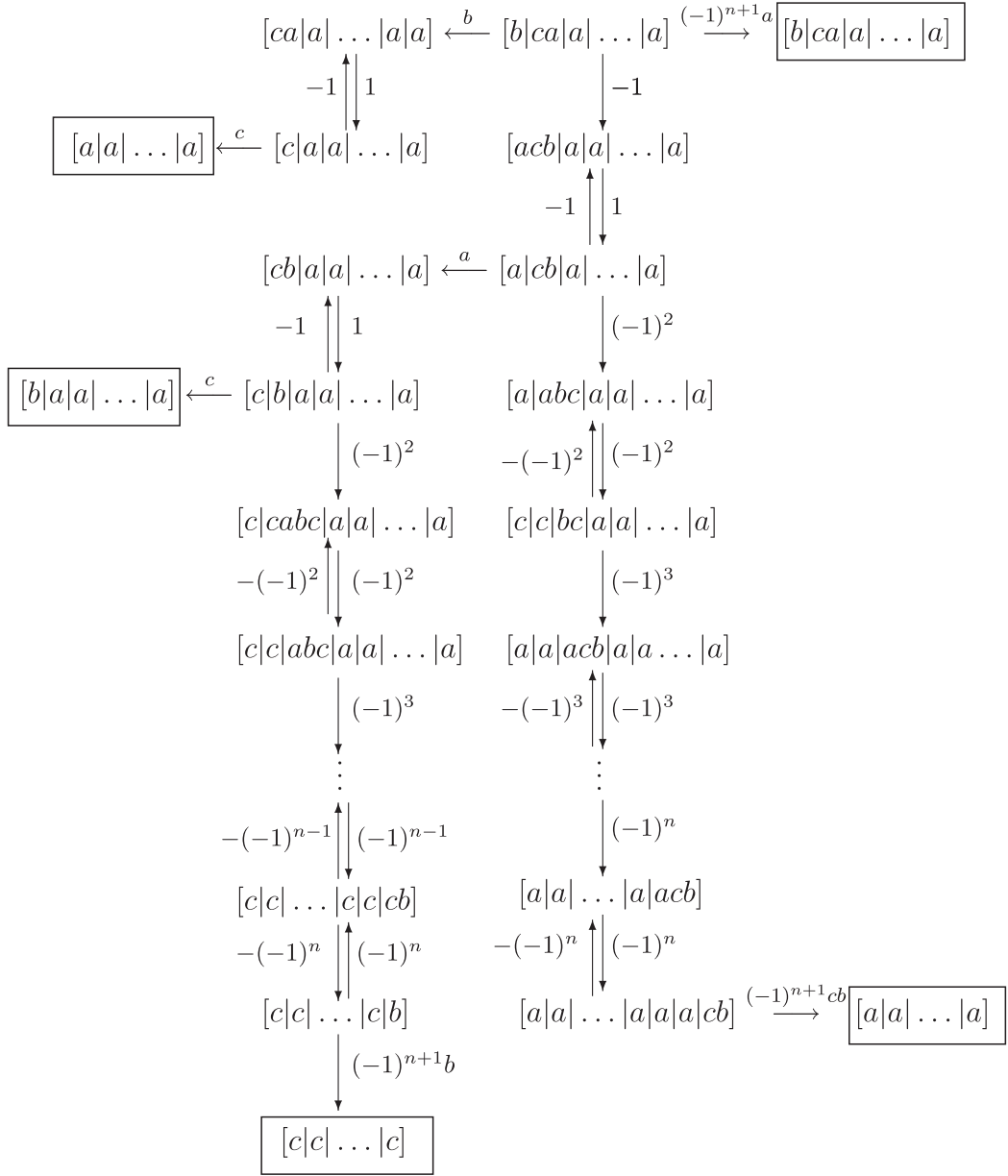}
\caption{Morse matching graph for $[bca^n]$}
\end{figure}

\begin{figure}\label{fig2}
\vspace{0.1in}
\includegraphics{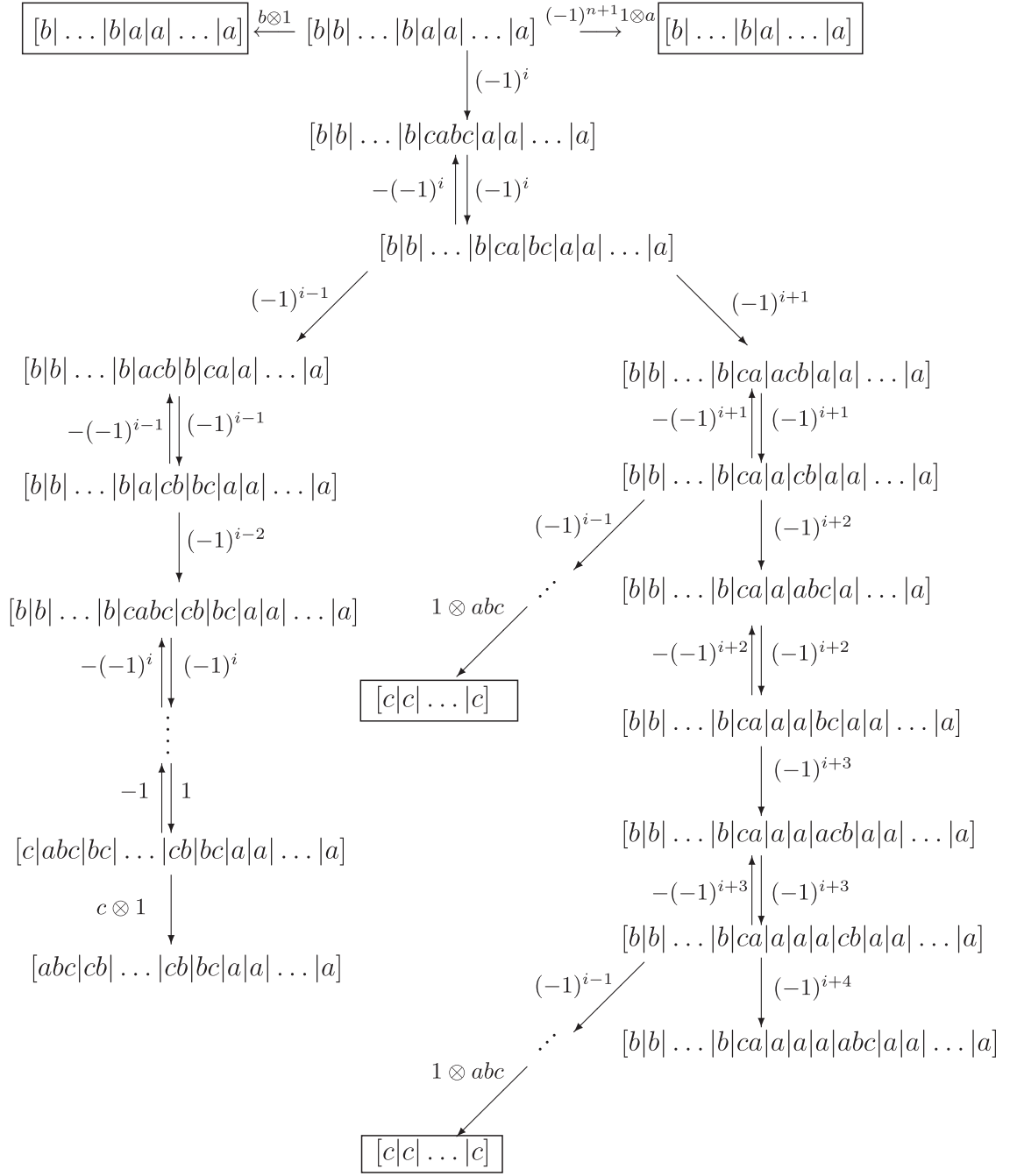}
\caption{Morse matching graph for $[b^ia^j]$, $n$ and $j$ are odd}
\vspace{1in}
\end{figure}

\begin{figure}\label{fig3}
\vspace{0.1in}
\includegraphics{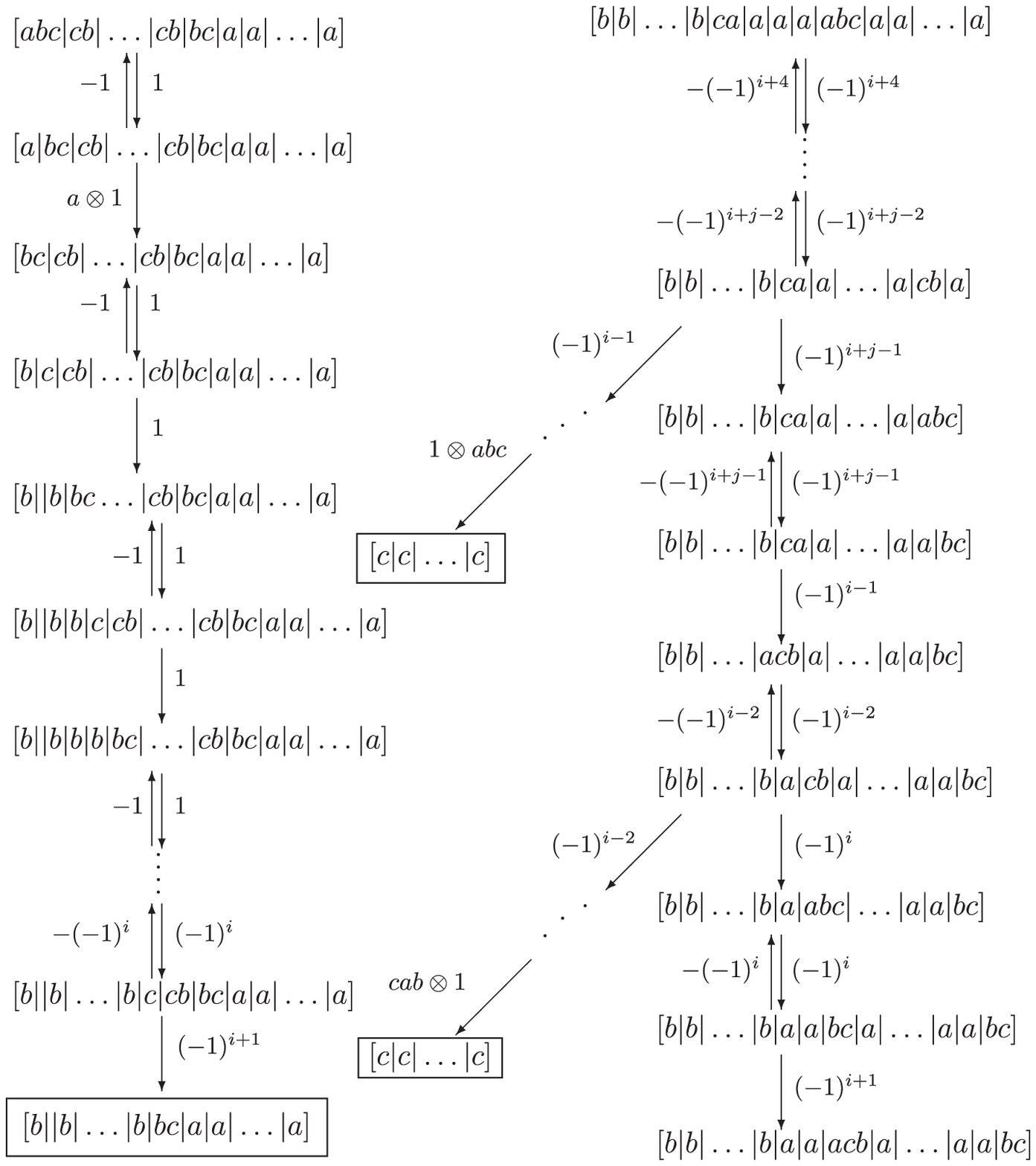}
\caption{Morse matching graph for $[b^ia^j]$, continued}
\vspace{1in}
\end{figure}

\begin{figure}\label{fig4}
\includegraphics{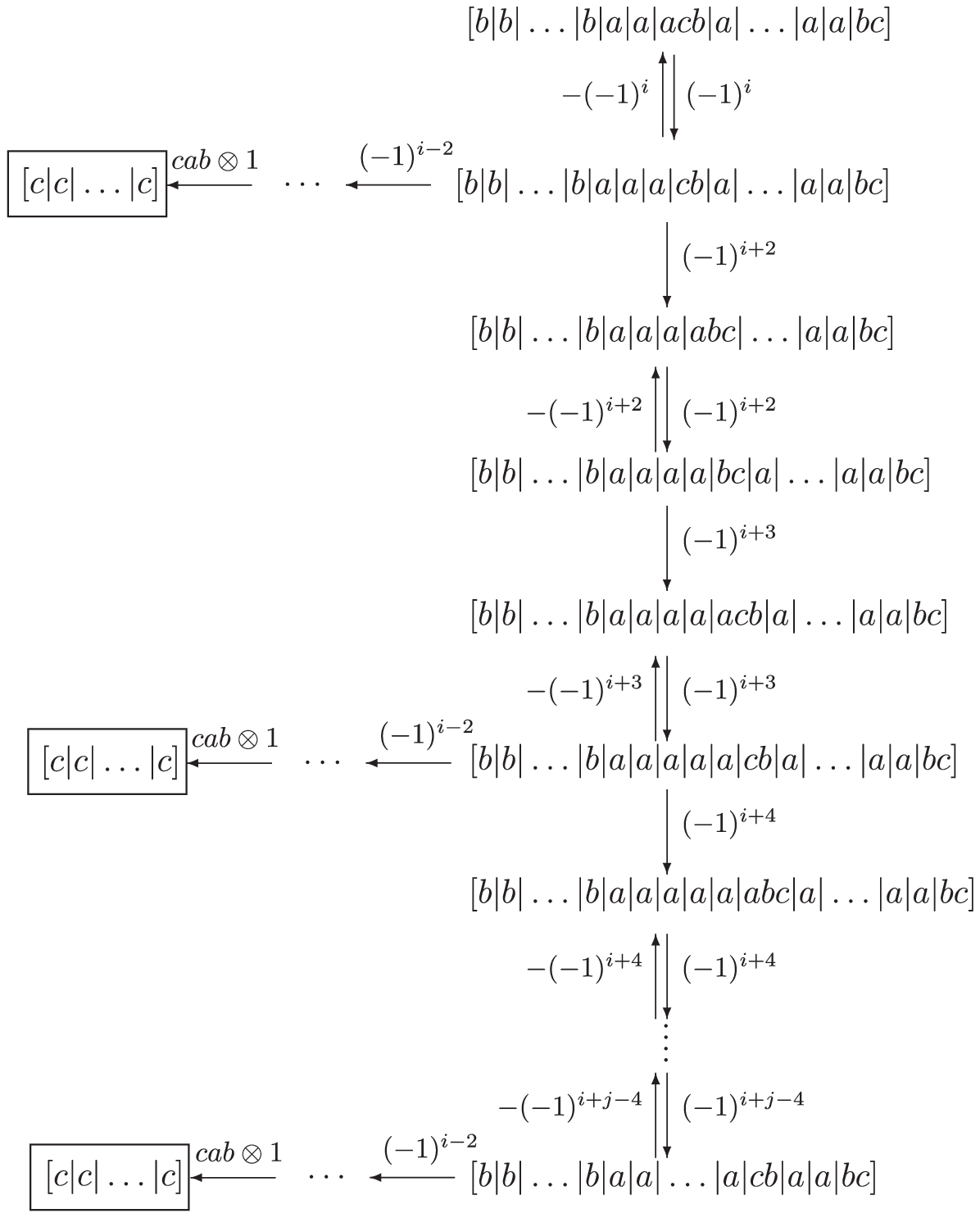}
\caption{Morse matching graph for $[b^ia^j]$, continued}
\end{figure}

Figure~\ref{fig1} 
presents the subgraph which allows us to calculate $d_{n+1}[bca^n]$.
For convenience, vertices from $V^{(n-1)}$ are boxed. 
For $[b^ia^j]$, $i+j=n+1$, $n+1$  is even, $j$ is odd, the corresponding 
graph is shown on Figures \ref{fig2}--\ref{fig4}. 
Other elements of $V^{(n)}$ can be processed in a similar way.\qed

\section{Hochschild cohomologies of $G_3^2$}

Suppose $W$ is a 1-dimensional bimodule over $kG^2_3$. 
Then $W=k$, and there exist 
$\epsilon_1,\epsilon_2: \{a,b,c\} \to \{1,-1\}$ such that 
$x1=\epsilon_1(x)$, $1x=\epsilon_2(x)$ for $x\in \{a,b,c\}$, $1\in k$.
Denote the bimodule obtained by
$W={}_{\epsilon_1}W_{\epsilon_2}$. 

It is natural to identify $\epsilon_i$ with a triple of signs, e.g., 
if $\epsilon_1(a)=\epsilon_1(b)=1$, $\epsilon_1(c)=-1$, 
$\epsilon_2(a)=1$, $\epsilon_2(b)=\epsilon_2(c)=-1$ then 
\[
 {}_{\epsilon_1}W_{\epsilon_2} = {}_{{+}{+}{-}}W_{{+}{-}{-}}.
\]

Our aim is to calculate $H^n(kG_3^2, {}_{\epsilon_1}W_{\epsilon_2})$ for all 
$\epsilon_i$, $i=1,2$. Up to the natural equivalence, it is enough to consider the following 
bimodules:
\begin{enumerate}
 \item 
$W_1 = {}_{{+}{+}{+}}W_{{+}{+}{+}}$;
 \item 
$W_2 = {}_{{+}{+}{+}}W_{{-}{-}{-}}$;
 \item 
$W_3 = {}_{{+}{+}{+}}W_{{+}{+}{-}}$; 
 \item 
$W_4 = {}_{{+}{+}{+}}W_{{+}{-}{-}}$; 
 \item 
$W_5 = {}_{{+}{+}{-}}W_{{+}{+}{-}}$; 
 \item 
$W_6 = {}_{{+}{+}{-}}W_{{-}{-}{+}}$; 
 \item 
$W_7 = {}_{{+}{+}{-}}W_{{+}{-}{-}}$; 
\item
$W_8 = {}_{{+}{+}{-}}W_{{+}{-}{+}}$.
\end{enumerate}

For every $W_i$ above, the map
 $d_{n+1}$ from Theorem~\ref{thm:Res} (for $n=2,3$, see Examples \ref{exmp:d_1}, \ref{exmp:d_2}) 
 induces a linear map 
from $kV^{(n)}$ to $kV^{(n-1)}$, let us denote it by $\hat d_{n+1}$. 
The corresponding conjugate map $\hat d_{n+1}^*$ acts from  
the space $(kV^{(n)})^*$ to $(kV^{(n+1)})^*$, and 
\[
H^{n+1}:=H^{n+1}(kG_3^2, W)\simeq \Ker \hat d_{n+2}^*/\mathrm{Im}\,\hat d_{n+1}^*, 
\]
so 
\begin{equation}\label{eq:HochDim}
\dim H^{n+1} = |V^{(n)}| - \rank \hat d_{n+2} -\rank \hat d_{n+1}.
\end{equation}
Since $|V^{(n)}|=2n+3$, it is enough to evaluate $\rank \hat d_{n+1}$ which is 
straightforward. 

The group $G_3^2$ is generated by involutions, so the first Hochschild cohomology group is 
expectedly trivial.

\begin{corollary}\label{cor1}
For $i=2,3,4,6,7,8$, we have $H^2(kG_3^2, W_i)=0$, but for
$i=1,5$ the second Hochschild cohomology group is nontrivial, 
$\dim H^2(kG_3^2, W_i)=1$.
\end{corollary}

\begin{corollary}\label{cor2}
For $n\ge 3$, $H^n(kG_3^2, W)=0$ for all 1-dimensional $kG_3^2$-bimodules.
\end{corollary}

\end{document}